\documentclass{amsart}
\usepackage{wha}
\usepackage{mathbbol}
\usepackage[all]{xy}
\author{Valery Alexeev}
\title{Weighted grassmannians and stable hyperplane arrangements}
\date{June 4, 2008}

\begin{document}

\begin{abstract}
  We give a common generalization of (1) Hassett's weighted stable
  curves, and (2) Hacking-Keel-Tevelev's stable hyperplane arrangements.
\end{abstract}
\maketitle

\tableofcontents

\section{Introduction and main statements}

The moduli space $\oM_{0,n}$ of stable $n$-pointed rational curves 
has many generalizations, beginning of course with $\oM_{g,n}$. For
this paper, however, the following two generalizations will be
important: 
\begin{enumerate}
\item Hassett's moduli $\oM_{0,\w}$ of weighted stable $n$-pointed
  curves \cite{Hassett_WeightedStableCurves}, and
\item Hacking-Keel-Tevelev's moduli $\oM(r,n)$ of stable hyperplane
  arrangements \cite{HackingKeelTevelev}.
\end{enumerate}

A \emph{weight data}, or simply a \emph{weight},
$\w$ is a collection of $n$ of rational (or real) numbers $0 <
b_i\le~1$. 
We denote $\bo=(1,\dotsc,1)$. 
A weighted stable curve of genus zero is a nodal curve
$X=\cup \bP^1$ whose dual graph is a tree, together with $n$ points
$B_1, \dotsc, B_n$ satisfying two conditions:
\begin{enumerate}
\item (on singularities) $B_i\ne$ the nodes, and whenever some points
  $\{B_i,\ i\in I\}$ coincide, one has $\sum_{i\in I}b_i\le~1$.
\item (numerical) $K_X+\sum b_i B_i$ is ample. In plain words, this
  means that for every irreducible component $E$ of $X$, one has
  \begin{math}
    |E\cap(X-E)| + \sum_{B_i\in E} b_i > 2.
  \end{math}
\end{enumerate}
The space $\oM_{0,\w}$ is the fine moduli space for flat families of
such weighted curves; it is smooth and projective.

A \emph{stable pair} $(X,B=\sum_{i=1}^n b_iB_i)$ is a natural
higher-dimensional analogue of the above notion. It consists of
a connected equidimensional projective variety $X$ together with $n$ Weil
divisors $B_i$ satisfying the following conditions
(see \cite{Alexeev_ICMtalk} for more details):
\begin{enumerate}
\item (on singularities) $X$ is reduced, and 
the pair $(X,B)$ is slc (semi log canonical),
  and
\item (numerical) $K_X+B$ is ample. 
\end{enumerate}

In \cite{HackingKeelTevelev} the authors construct a projective
scheme, which we will denote $\oM(r,n)$, together with a flat family
$f:(\cX,\cB_1,\dotsc,\cB_n)\to\oM(r,n)$ such that every geometric
fiber $(X,\sum B_i)$ is a stable pair in the above definition, with
all coefficients $b_i=1$.  Over an open (but not dense in general) subset
$\M(r,n) \subset \oM(r,n)$ this gives the universal family of $n$
hyperplanes $B_i$ on a projective space $X=\bP^{r-1}$ such that $B_i$
are in general position.
The construction originates in
\cite{Kapranov_QuotientsGrassmannians}, see also \cite{Lafforgue_Chirurgie}.

More generally, let $\w$ be a weight, and $B_1,\dotsc, B_n$ be $n$
hyperplanes in $\bP^{r-1}$. Then the pair $(\bP^{r-1}, \sum b_i B_i)$
is 
\begin{enumerate}
\item \emph{lc (log canonical)} if for each intersection $\cap_{i\in I}B_i$
  of codimension $k$, one has $\sum_{i\in I}b_i \le k$, and
\item \emph{klt (Kawamata log terminal)} if the inequalities are
  strict, in particular all $b_i<1$.
\end{enumerate}
(This is consistent with the standard definitions of the Minimal Model
Program.) 

The pair $(\bP^{r-1}, \sum b_i B_i)$ is stable in the above definition
iff it is lc (slc being an analog of lc for possibly nonnormal pairs) and
$|\w|= \sum_{i=1}^n b_i > r$.  We call such pairs \emph{weighted
  hyperplane arrangements}, or simply \emph{lc hyperplane
  arrangements}. 
One easily constructs a fine moduli space
$\M_{\w}(r,n)$ for them; it is smooth, of dimension $(r-1)(n-r-1)$,
and usually not complete (but see Theorem \ref{thm:small-weights} for
the exceptions).

Throughout the paper, we work over an arbitrary commutative ring $\cA$
with identity. The main results of this paper are the three theorems
below and the detailed description of the weighted stable hyperplane
arrangements given in Section~\ref{sec:defin-moduli-univ}.

\begin{theorem}[Existence]\label{thm:existence}
  For each $r,n$ and a rational weight $\w=(b_i)$ with $|\w|=\sum
  b_i>r$, there exists a projective scheme $\oM_{\w}(r,n)$ together
  with a locally free (in particular, flat) family
  $f:(\cX,\cB_1,\dotsc,\cB_n) \to \oM\lw(r,n)$ such that:
  \begin{enumerate}
  \item Every geometric fiber of $f$ is an $(r-1)$-dimensional variety
    $X$ together with $n$ Weil divisors $B_i$ such that the pair
    $(X,\sum b_i B_i)$ is stable.
  \item For distinct geometric points of $\oM\lw(r,n)$, the fibers are
    non-isomorphic. 
  \item Over an open (but not dense in general) subset
    $\M_{\w}(r,n)\subset \oM_{\w}(r,n)$, $f$ coincides with the universal
    family of weighted hyperplane arrangements.
  \end{enumerate}
  For every positive integer $m$ such that all $mb_i\in\bN$, the sheaf
  $\cO_{\cX}(m( K_{\cX}+\sum b_i\cB_i) )$ is relatively ample and
  free over $\oM\lw$.
\end{theorem}

The fibers of $f$ will be called \emph{weighted stable hyperplane
  arrangements}, or simply \emph{slc hyperplane arrangements}.  
As one has $\M_{\bo}(r,n) \subset
\M_{\w}(r,n)$, in particular, each of $\oM_{\w}(r,n)$ provides a
moduli compactification of the moduli space of generic hyperplane
arrangements.

\begin{definition}\label{defn:chambers}
We define the \emph{weight domain} 
\begin{displaymath}
  \cD(r,n) = \{ (b_i)\in \bQ^n \mid 
  0<b_i\le 1,\quad \sum b_i > r \}
\end{displaymath}
and a subdivision of it into locally closed \emph{chambers},
denoted $\chamber(\w)$, by the
hyperplanes $\sum_{i\in I}b_i = k$ for all $1\le k \le r-1$ and
$I\subset \{1,\dotsc,n\}$, and by the faces $b_i=1$.

We introduce a partial order on the points of $\cD(r,n)$: $\w>\w'$ if
for all $1\le i \le n$ one has $b_i\ge b_i'$, with at least one
strict inequality. 
\end{definition}

We will frequently assume $(r,n)$ fixed and drop it from the notation.

\begin{example}\label{exmp:3-degs}
  ($r=3$, $n=5$) Consider a 1-parameter family of 5 lines on $\bP^2$
  in general position such that in the limit $B_1,B_2,B_5$ meet at a
  point $q_1$, and $B_3,B_4,B_5$ meet at a point $q_2\ne q_1$. This is
  not allowed by the lc
  singularity condition if $b_1+b_2+b_5>2$ or 
  $b_3+b_4+b_5>2$. Since the spaces $\oM_{\w}$ are proper, there is
  always a stable pair limit, but its shape depends on the weight:
  \begin{enumerate}
  \item $\w=(1,1,1,1,1-\epsilon)$. The variety is $X_{\w}=X^0\cup X^1\cup
    X^2$, where $X^0$ is the blowup of $\bP^2$ at 
    $q_1$ and $q_2$, and $X^0$ is glued along the exceptional $\bP^1$'s to
    $X^1=\bP^2$, $X^2=\bP^2$.  The divisor $B_5$ has three
    irreducible components, each of $B_1,\dotsc, B_4$ has two.
    $B_1,B_2$ are contained in $X^0\cup X^1$, and $B_3,B_4$ in
    $X^0\cup X^2$. All five divisors are Cartier.
  \item $\w'=\bo=(1,\dotsc,1)$. The variety $X_{\w'}$ is obtained from
    $X_{\w}$ by contracting the $(-1)$-curve $B_5\cap X^0$. The image
    of $X^0$ is ${X'}^0=\bP^1\times\bP^1$, the divisors $B_1,B_2$
    restricted to ${X'}^0$ are fibers of a ruling, and $B_3,B_4$
    restricted to it are fibers of the second ruling. The divisor
    $B_5$ intersects ${X'}^0$ at one point, and so is not $\bQ$-Cartier.
  \item $\w''=((1+\epsilon)/2,(1+\epsilon)/2,1,1,1-\epsilon)$. The
    variety $X_{\w''}$ is obtained from $X_{\w}$ by contracting $X^1$.
    \end{enumerate}
    Note that $\w'>\w>\w''$ and $\w',\w''\in \overline{\cham(\w)}$, 
    we have natural morphisms
    $X_{\w'}\gets X_{\w} \to X_{\w''}$, and the first of these
    morphisms is birational.
\end{example}

\begin{theorem}[Reduction morphisms]
\label{thm:red-moduli}
  \begin{enumerate}
  \item (Same chamber) For $\w,\w'$ lying in the same chamber, one has
    $\oM_{\w}=\oM_{\w'}$ and
    $(\cX,\cB_i)_{\w}=(\cX_,\cB_i)_{\w'}$. In particular, the divisors
    $\sum (b_i-b_i)B_i$ are $\bQ$-Cartier.
  \item For $\w'\in \overline{\chamber (\w)}$, there are
    natural reduction morphisms
    \begin{equation*}
      \xymatrix{
        (\cX,\cB_i)_{\w} \ar[r]^{\pi_{\w,\w'}} \ar[d] &
        (\cX,\cB_i)_{\w'} \ar[d] \\
        \oM_{\w} \ar[r]^{\rho_{\w,\w'}} &
        \oM_{\w'}
      }
    \end{equation*}
    One has
    \begin{displaymath}
      \pi_{\w,\w'}^* \cO_{\cX\lwp}(m( K_{\cX\lwp}+\sum b_i'\cB_i) ) 
      = \cO_{\cX\lw}(m( K_{\cX\lw}+\sum b_i'\cB'_i) )
    \end{displaymath}
    for any $m$ such that all $mb_i'\in\bZ$.
  \item (Specializing up) For $\w' \in \overline{\chamber (\w)}$ with
    $\w<\w'$, $\rho_{\w,\w'}$ is an isomorphism, and on the fibers
    $\pi_{\w,\w'}:X\to X'$ is a birational contraction restricting to
    an isomorphism $X\setminus\cup B_i \isoto X'\setminus\cup B'_i$.
  \item For any $\w>\w'$, there is a natural reduction morphism
    $    \rho_{\w,\w'}:\oM_{\w}\to \oM_{\w'}.    $
    On the fibers, the rational map $\pi_{\w,\w'}:X\dashrightarrow X'$ is a
    sequence of log crepant contractions and log crepant birational
    extractions. Further, $X'=\Proj \oplus_{d\ge0}H^0(X,\cO(dm(K_X+\sum
    b_i'B_i)))$ is the log canonical model for the pair $(X,\sum b_i'B_i)$.
  \end{enumerate}
\end{theorem}

\begin{theorem}[Moduli for small weights]
  \label{thm:small-weights}
  Let $\alpha=(a_i)$ be a weight with $\sum a_i=r$ (lying on the boundary of        
  $\cD$) which belongs to the closure of a unique chamber $\cham(\w)$. Then
  $$\M\lw=\oM_{\w} =(\bP^{r-1})^n//PGL(r) = \Gr(r,n) //\bG_m^{n-1} $$
  is the GIT quotient for the line bundle, resp.
  linearization corresponding to $\alpha$.
\end{theorem}

For \emph{any} boundary weight (i.e. with $|\alpha|=r$), we can formally
define $\oM_{\alpha}$ to be the above GIT quotient.  Over an open
and dense subset $\M_{\alpha}$ it gives the moduli of lc hyperplane
arrangements on $\bP^{r-1}$ such that $K_{\bP^{r-1}} + \sum
a_i B_i=0$. For $\alpha$ as in the theorem, one has $\M_{\alpha}=\oM_{\alpha}
=\oM\lw$.


\begin{notations}
  We work over an arbitrary commutative base ring $\cA$ with identity,
  without the Noetherian assumption, and indeed can work over any base
  scheme. $A$ will denote an $\cA$-algebra, and $k=\bar k$ an
  $\cA$-algebra which is an algebraically closed field.  The tilde
  will be used to denote affine schemes $\wt X$, cones $\wt\Delta$,
  etc., which are cones over the corresponding projective schemes $X$,
  polytopes $\Delta$, etc.
\end{notations}

It may help the reader to grasp some combinatorial aspects of this
paper with the following general outline. The (unweighted) stable
hyperplane arrangements are described by matroid tilings of the
hypersimplex $\Delta(r,n)$. Their weighted counterparts are described
by partial tilings of $\Delta(r,n)$ as viewed through a smaller
``window'' $\Delta_{\w}(r,n)$; the window must be completely covered.

Another key idea is the GIT interpretation of the weight $\beta$
explained in Section~\ref{sec:git}.

\section{Matroid polytopes}\label{sec:matroid-polytopes}

We begin with some general definitions and then specialize them to the
case of grassmannians.

\begin{setup}\label{setup}
  We fix two lattices $\bZ^N=\oplus\bZ e_j$ and $\bZ^n$, a homomorphism
  $\phi:\bZ^N\to \bZ^n$, and a homomorphism $\deg:\bZ^n\to\bZ$, such
  that $\deg \phi(e_j)=1$ for all $j$.
  Associated to this data are affine $\bA=\bA^N$ and
  projective $\bP=\bP^{N-1}$ spaces over $\cA$
  and linear actions of split tori $\wT=\bG_m^n$ on $\bA$ and of
  $T=\wT/\diag\bG_m$ on $\bP$.

  Let $\Delta$ be the lattice polytope that is the convex hull of
  $\phi(e_i)$, and $\wt\Delta$ be the corresponding cone in $\bR^n$.
  We also fix a $\bZ^n$-graded ideal
  $I[\wZ]\subset\cA[z_1,\dotsc,z_N]$ such that the quotient is a locally
  free (i.e. projective) $\cA$-module.
  Hence, $\wZ\subset\bA$ is a $\wT$-invariant closed subscheme. Let
  $Z\subset \bP$ be the corresponding $T$-invariant closed subscheme.
\end{setup}

\begin{definition}
For a geometric point $p\in Z(k)$, the closure of the orbit
$\overline{T.p}$ is a possibly nonnormal toric subvariety of
$Z_k=Z\times_{\cA}k$. 
It corresponds to a
lattice polytope $P$ which we will call the \emph{$Z$-polytope} or the
\emph{moment polytope of $p$}. 
(Indeed, when $k=\bC$, $P$ is  the moment polytope of
$\overline{T.p}$, as defined in symplectic geometry.) A character
$\chi\in\bZ^n$ is in the cone $\wP$ iff there exists a monomial
$z^m=\prod_i^N z_i^{m_i}$ such that $\phi(m)=d\chi$ and $z^m(p)\ne0$.

A \emph{$Z$-tiling} $\uP$ is a face-fitting subdivision of $\Delta$
into $Z$-polytopes.
\end{definition}

We fix several faces $F_i$, $i=1,\dotsc,n'$, of $\Delta$. Each of them
is defined by the inequality $l_i\le1$ for a unique $\bZ$-primitive
linear function $l_i(x_1,\dotsc,x_n)$. In a completely parallel
fashion with our grassmannian setup, an element $\w=(b_i)\in\bQ^{n'}$,
$b_i\le1$, is called a weight.
For each weight we define a subpolytope 
\begin{displaymath}
  \Delta\supset\Delta\lw = \{ l_i\le b_i \}
\end{displaymath}
The \emph{weight domain} $\cD\subset\bR^{n'}$
is the set of the weights for which $\Delta\lw$ is nonempty and
maximal-dimensional. 

\begin{definition}
  A \emph{weighted $Z$-polytope} is a polytope of the form $P\lw=
  P\cap \Delta_{\w}$ for some $Z$-polytope $P$, called 
  \emph{the parent} of $P\lw$, such that $\Int(P)\cap
  \Delta_{\w}\ne \emptyset$.

  A \emph{weighted $Z$-tiling} $\uP\lw$ is a face-fitting
  tiling of $\Delta_{\w}(r,n)$ by weighed $Z$-polytopes.
  The \emph{partial} cover $\uP$ of $\Delta$ by the parent
  polytopes is called \emph{the parent cover} of $\uP\lw$.  
\end{definition}

\begin{definition}\label{defn:Z-chamber-decomp}
  The $Z$-chamber decomposition of $\cD$ is defined as follows: $\w,\w'$
  lie in the same chamber if for every $Z$-polytope $P$, one has
  $P\cap\Delta\lw\ne\emptyset \iff P\cap\Delta\lwp\ne\emptyset$.
  Consequently, weighted $Z$-tilings of $\Delta\lw$ and $\Delta\lwp$
  are in a bijection.
\end{definition}

We now specialize these definitions to the case of the
grassmannians. The polytope $\Delta$ in this case is called the
\emph{hypersimplex} and the $Z$-polytopes are called \emph{matroid
  polytopes}. For the unweighted version, these notions were
introduced in \cite{GelfandGoreskyMacPhersonSerganova}.

Let $\Gr(r,n)$ be the grassmannian of $r$-planes in a fixed affine
space $\bA^n$, together with its Pl\"ucker embedding into
$\bP(\wedge^r \bA^n)=\bP^{N-1}$, where   $N={n\choose r}$.
Let $\wgr(r,n)\subset \bA^N$  be the affine cone. It is defined by the
classical quadratic Pl\"ucker relations.

For $I=(i_1,\dotsc,i_r)$, the Plu\"cker coordinate $p_I$ has character
$$\weight(p_I)=(1,\dotsc,1,0\dotsc,0)$$ with $r$ ones in the places
$i_1,\dotsc,i_r$ and with $(n-r)$ zeros elsewhere.

\begin{definition}
  The convex hull of these $N$ points is called the
  \emph{hypersimplex} $\Delta(r,n)$. Alternatively, it can be
  described as follows:
  \begin{displaymath}
    \Delta(r,n) = \left\{ (x_i)\in \bR^n \mid 0\le x_i \le 1, \quad 
    \sum_{i=1}^n x_i = r \right\}
  \end{displaymath}
It has $2n$
faces $F_i^+= \{x_i=1\}$ and $F_i^-=\{x_i=0\}$, isomorphic to
$\Delta(r-1,n)$ and $\Delta(r,n-1)$ respectively.

We fix the lattice $\Lambda\simeq\bZ^n$ in $\bR^n$ consisting
of the vectors $(x_i)\in\bZ^n$ such that $\sum x_i$ is divisible by
$r$ and a homomorphism $\deg:\Lambda\to\bZ$ so that the characters of the
Pl\"ucker coordinates $p_I$ generate $\Lambda$ and have degree 1.
\end{definition}

\begin{definition}
  A \emph{matroid polytope} $P_V\subset \Delta(r,n)$ is the polytope
  corresponding to the toric variety $\overline{T.V}$ for some
  geometric point $[V\subset\bA^n]\in\Gr(r,n)(k)$.
  (Theorem~\ref{thm:matroid-polytopes}(1) implies that this projective
  toric variety and the corresponding affine variety are normal,
  unlike the case of general $Z$).
\end{definition}

The equations of the coordinate hyperplanes restricted to $V$ give
$n$ vectors $z_1, \dotsc, z_n \in V^*$, what is called a
\emph{realizable matroid}.
Then $\weight(p_{i_1,\dotsc,i_r})$ is a vertex of $P_V$ iff
$z_{i_1}, \dotsc, z_{i_r}$ form a basis of $V^*$. Alternatively,
$P_V$ can be described inside $\Delta(r,n)$ by the inequalities
$\sum_{i\in I}x_i \le \dim\ \Span(z_i,\ i\in I)$, for all $I\subset
\{1,\dotsc,n\}$.

One can also describe the matroid polytopes in terms of hyperplane
arrangements. Let $\bP V\simeq\bP^{r-1}$ be the corresponding
projective space and assume that it is not contained in any of the $n$
coordinate hyperplanes $H_i$ (i.e. all $z_i\ne0$ on $\bP V$); let 
$B_1,\dotsc,B_n\subset \bP V$ be $H_i\cap\bP V$. Then 
$\weight(p_{i_1,\dotsc, i_r})$ is a vertex of $P_V$ iff
$B_{i_1}\cap\dotsc \cap B_{i_r}$ is a point. Alternatively, 
$P_V$ can be
described inside $\Delta(r,n)$ by the inequalities $\sum_{i\in I}x_i
\le \codim\cap_{i\in I}B_i$ for all $I\subset \{1,\dotsc,n\}$. Note
that the matroid polytope in this case is not contained in any of the faces
$F_i^-=\{x_i=0\}$. 

\begin{definition}
  A \emph{matroid tiling} $\uP$ of $\Delta(r,n)$ is a face-fitting
  subdivision $\cup P(V_{s})$  of $\Delta(r,n)$ into matroid polytopes.
\end{definition}

Matroid polytopes form a very particular class of lattice polytopes,
with many properties not shared by general lattice polytopes. Some of their
properties can be summarized as follows:

\begin{theorem}\label{thm:matroid-polytopes}
  \begin{enumerate}
  \item Every matroid polytope is generating, i.e. its integral points
    generate the group of integral points of $\bR P$. Moreover, the
    semigroup of integral points in $\wP$ is generated by the vertices
    of $P$.
  \item A matroid polytope of codimension $c$ is in a canonical way
    the product of $c+1$ maximal-dimensional matroid polytopes for smaller
    $r,n$. So, one has $r=r_0+\dotsc+ r_c$ and 
    $\{1,\dotsc,n\}=I_0\sqcup \dotsb \sqcup I_c$, and $P=\prod P_j$,
    where $P_j\subset \Delta(r_j,|I_j|)$ is a maximal-dimensional
    matroid polytope.
  \end{enumerate}
\end{theorem}

We now introduce the weighted versions of these notions. 

\begin{definition}
  Let $\w=(b_1,\dotsc,b_n)$ be a weight. A \emph{weighted
    hypersimplex} is the polytope given by 
  \begin{displaymath}
    \Delta_{\w}(r,n) = \left\{ (x_i)\in \bR^n \mid 0\le x_i \le b_i, \quad 
    \sum_{i=1}^n x_i = r \right\}.
  \end{displaymath}
  Similarly, we also have definitions of a \emph{weighted matroid
  polytope}, a \emph{weighted matroid tiling $\uP\lw$}, and \emph{the
  parent cover} of $\uP\lw$.
\end{definition}

\begin{question}\label{qu:extend-to-parent}
  Can every parent cover be extended to a \emph{complete cover} of
  $\Delta(r,n)$? For $r=2$ the answer is easily seen to be ``yes''.
  For $r\ge3$ we expect the answer to be ``no'', following the general
  philosophy of ``Mnev's universality theorem'' (cf.
  \cite[Thm.I.14]{Lafforgue_Chirurgie} which shows that matroid geometry
  can be arbitrarily complicated.
\end{question}

\begin{theorem}[Chamber decomposition]
  \label{thm:chamber-decs-coincide}
  The chamber decomposition of $\cD(r,n)$ defined in
  \ref{defn:Z-chamber-decomp} coincides with that of
  Definition \ref{defn:chambers}.
\end{theorem}
\begin{proof}
  Starting with $\w=\w'$ and then varying the weight $\w$, the matroid
  decompositions of $\Delta\lw$ and $\Delta\lwp$ may possibly change if for
  some matroid polytope $P\subset\Delta(r,n)$, a vertex of $\Delta\lw$
  would lie on a face of $P$ that $\w'$ did not belong to.

  Every vertex of $\Delta\lw$ is of the following form: $x_i=0$ or
  $b_i$, for all but possibly one ${i_0}$.  Every face of a matroid
  polytope lies in the intersection of hyperplanes $\sum_{i\in I}x_i =
  k$ for some $1\le k\le n-1$ and $I\subset\{1,\dotsc,n\}$. 
  Possibly after replacing
  $I$ by its complement, we can assume that $i_0\not\in I$. Then for
  some $J\subset I$ we get $\sum_{i\in J}b_i=k$. If $\w$ belongs to a
  face of $P$ that $\w'$ did not belong to, then we get a new equation
  of this form. So $\w$ lies in a different, smaller chamber.
\end{proof}

\cite{GelfandGoreskyMacPhersonSerganova} gives three different
interpretations of matroid polytopes. Here, we add another one.

\begin{theorem}\label{thm:matroid-lc}
  The matroid polytope $P_V$ is the set of points $(x_i)\in \bR^n$
  such that the pair $(\bP V,\sum x_iB_i)$ is lc and 
  $K_{\bP V}+\sum x_iB_i=0$; 
  the interior $\Int P_V$ is the set of points such that
  $(\bP V,\sum x_iB_i)$ is klt and   $K_{\bP V}+\sum x_iB_i=0$.
\end{theorem}
\begin{proof}
  Indeed, the defining inequalities $\sum_{i\in I}x_i\le \codim \cap_{i\in
    I}B_i$ of $P_V$ also happen to be the conditions for the   pair $(\bP
  V,\sum x_iB_i)$ to be lc. Similarly with the strict inequalities and klt. 
\end{proof}

We don't even have to assume that $\bP V\not\subset H_i=\{z_{i_0}=0\}$
for this theorem: the pair $(\bP V,\sum x_iB_i)$ can only be lc if
$x_{i_0}=0$, otherwise $\sum x_iB_i$ is not a divisor. And indeed if
$\bP V\subset H_i$ then $P_V\subset \{x_i=0\}$, so the theorem
still holds.

\section{Moduli spaces for varieties with torus action}
\label{sec:two-moduli-spaces}

Let $Z\subset \bP$ be a projective scheme locally free over $\cA$ and
invariant under the $T$-action, and $\wt Z\subset\bA$ be its affine cone,
with the $\wT$-action, as in our general setup \ref{setup}. 
Two moduli spaces of varieties with torus
action will be relevant for this paper. 

\begin{enumerate}
\item The toric Hilbert scheme $\Hilb^{\wT}(\wZ,\wt\Delta)$, constructed
  in \cite{PeevaStillman,HaimanSturmfels}.
\item The moduli space $\M^T(Z,\Delta)$ of finite $T$-equivariant maps
  $Y\to Z$ of stable toric varieties $Y$ over $Z$,
  constructed in \cite{AlexeevBrion_SphericalModuli}, see also
  \cite{AlexeevBrion_Affine,AlexeevBrion_Projective,Alexeev_CMAV}.
  This is the equivariant multiplicity-free version of the moduli
  space of branchvarieties   \cite{AlexeevKnutson}.
\end{enumerate}
Both of these moduli spaces are projective schemes. Both
are available in much more general
settings; we will only need the simplest versions.

For an $\cA$-algebra $A$, $\Hilb^{\wT}(\wZ,\wt\Delta)(A)$ is the set
of closed $\wT$-invariant subschemes $\wt Y \subset
\wZ_A=\wZ\times_{\cA} A$ which are
\emph{multiplicity-free:} for
every $x\in\bZ^n$ the graded piece $A[\wY]_x$ is a locally free rank-1
$A$-module if $x\in\wt\Delta$ and is $0$ otherwise. A geometric fiber
$\wY_k$ need not be reduced. $(\wY_k)_{\rm red}$ is a
union of possibly non-normal toric varieties glued along torus orbits 
in a fairly complicated way.

In contrast, $\M^T(Z,\Delta)(A)$ is the set of locally free proper families
$Y$ over $\Spec A$ together with a finite $T$-equivariant morphism
$f:Y\to Z_A$ such that every geometric fiber $Y_k$ is a projective
\emph{stable toric variety.} In the ring $\oplus_{d\ge0}
H^0(Y,f^*\cO(d))$, for each $x\in\bZ^n$, the $x$-graded piece is a
locally free rank-1 $A$-module if $x\in\wt\Delta$ and is $0$
otherwise.  A stable toric variety is a \emph{reduced} variety glued
from \emph{normal} toric varieties along orbits in a very simple way,
so that the result is seminormal.  The price for such niceness is that $Y_k\to
Z_k$ is a finite morphism rather than a closed embedding.

The projective stable toric variety $Y_k$ comes with the polarization
$L=f^*\cO_{Z}(1)$. For each irreducible component of $Y_k$, this
gives a lattice polytope $P^{s}$, and together they give a
tiling $\Delta=\cup_{s\in \uP} P^{s}$ describing the gluing
in a rather precise way. 
If the cone semigroups $\wP^{s}\cap\bZ^n$ are generated in degree
1 then $Y_k\to Z_k$ is a closed embedding and gives a point of
$\Hilb^{\wT}(\wZ,\wt\Delta)$.

In general, $\M^T(Z,\Delta)$ is only a coarse moduli space since a
finite map $Y_k\to Z_k$ may have automorphisms (deck
transformations). However, it is a fine moduli space over an open
subscheme where $Y_k\to Z_k$ is birational to its image (on every
irreducible component).

Finally, we note how these moduli spaces change if we replace $\bP$ by
a $d$-tuple Veronese embedding, which means replacing the ring
$\cA[z_1,\dotsc,z_N]$ by the subring generated by monomials of degrees
divisible by $d$. The answer for $\Hilb^{\wT}(\wZ)$ is very
non-obvious, and sometimes they indeed change. The moduli space
$\M^T(Z)$, however, does not change. Indeed, the scheme $Z$ does not
change, and neither does the $T$-action.  Therefore, $\M^T(Z)$ can be
defined as easily for a \emph{rational polytope} $\Delta$: it can
always be rescaled to make it integral.

\section{Review of the unweighted case}
\label{sec:unweighted-case}

We now apply the general theory of the previous section to the
grassmannians. Let $\gr=\gr(r,n)$ be the grassmannian with its
Pl\"ucker embedding into $\bP^N$, $N={n \choose r}$, and $\wgr$ be its
affine cone. Hence, $\cA[\wgr]$ is generated by the $N$ Pl\"ucker
coordinates $p_{i_1,\dotsc,i_r}$, modulo the usual quadratic
relations. The corresponding polytope is precisely the hypersimplex
$\Delta(r,n)$, and the polytopes $P^{s}$ appearing in the
constructions of the previous section are the matroid polytopes.

\begin{definition}[\cite{HackingKeelTevelev}]
  $\oM(r,n) = \Hilb^{\wt  T}(\wgr,\wt\Delta).$
\end{definition}

A closed $\wT$-invariant multiplicity-free subscheme $Y_k\subset
\gr_k$ gives a matroid subdivision of $\Delta(r,n)$.
\ref{thm:matroid-polytopes}(1) implies that $Y_k$ is reduced and
is a stable toric variety, so we are in the ``nice case''.  Note that $\dim
Y=n-1\ne r-1$, so $Y$ is not the required stable variety $X$. Instead,
it should be thought of as the \emph{log Albanese variety} $\logAlb(X,B)$.

\begin{definition}[\cite{HackingKeelTevelev}]
  $X\subset Y$ is the intersection of $Y$ with the subvariety
  in $\Gr(r,n)$ defined by
  \begin{displaymath}
    \gr^e = \{ V\subset \bA^n \mid (1,1,\dotsc,1) \in V \}.
  \end{displaymath}
  $\gr^e$ is a Schubert variety, isomorphic
  to $\Gr(r-1,n-1)$. 
\end{definition}

One easily shows that $\gr^e \into \Gr$ is a regular codimension $n-r$
embedding (not $\bG_m^n$-equivariant), the zero set of a section of the
tautological bundle $Q$ on $\Gr$.  $X$ does not contain any
$T$-orbits. This implies that $X\subset Y$ is a regular codimension
$n-r$ embedding as well, the zero section of the bundle $Q|_Y$ with
$c_1(Q|_Y)=L$.

On the other hand, by \cite{Alexeev_CMAV}, the pair $(Y,\sum
B_i^{\pm})$ has semi log canonical singularities and $K_Y+\sum
B_i^{\pm}=0$, where $B_i^{\pm}$ are the divisors corresponding the
faces $F_i^{\pm}$ of $\Delta(r,n)$. In addition, $B_i^-\cap
X=\emptyset$. Denoting $B_i^+|_X= B_i$ and combining this together
gives the following:

\begin{theorem}[\cite{HackingKeelTevelev}]
  \label{thm:HKT-main}
  \begin{enumerate}
  \item There exists a smooth morphism $X\times T \to Y$ whose image
    is the open subset $Y\setminus \cup B_i^-$ swept by the $T$-orbits
    of $X\subset Y$; this is compatible with the divisors $B_i$ and $B_i^+$.
  \item Consequently, $(X,\sum B_i)$ and $(Y\setminus \cup B_i^-, \sum
    B_i^+)$ are isomorphic locally in smooth topology; in particular,
    the pair $(X,\sum B_i)$ is slc.
  \item $K_X + \sum B_i = (K_Y + \sum B_i^{\pm} + L)|_X = L|_X$; and
    so is ample.
  \item The poset of the stratification of $X$ defined by the
    irreducible components, their intersections, and the divisors
    $B_i$ coincides with the poset of the stratification on
    $\Delta(r,n)\setminus \cup F_i^-$ defined by the subdivision $\uP$.
  \end{enumerate}
\end{theorem}

\section{Weighted grassmannians}
\label{sec:weight-grassm}

Here we define certain projective schemes $\gr\lw$ and describe their
basic properties. We begin with the elementary case which already
contains the pertinent combinatorics of the general situation.

Let $P'$ be a lattice polytope, and $Y$ be the correponding projective
toric scheme over $\cA$ (a toric variety when working over a field $k$),
together with an ample $T$-linearized ample invertible sheaf $L'$.
Let $m$ be a positive integer, $P=P'/m$ a rational polytope, and 
$L=L'/m \in \Pic(Y)\otimes\bQ$
be the corresponding $\bQ$-polarization.

Let us fix several faces $F_i$ of this polytope. Each of
them is given by a linear equation $x_i=b_i\in \bQ$, where $x_i$ is an
integral primitive linear function, so that $P$ lies in the half space
$x_i\le b_i$. Each of these faces corresponds to a 
divisor $B_i$ on $Y$. 

Now consider the polytope $P_{\w'}$ obtained by replacing the inequalities
$x_i\le b_i$ with  $x_i\le b_i'$ for some $b_i'\in\bQ$. 
Here, we denote $\w=(b_i)$ and $\w'=(b_i')$, so that $P=P_{\w}$.
Note that for some $\w'$ one may have $\dim P_{\w'}< \dim P_{\w}$. 

One says that
two polytopes are \emph{normally equivalent} if their normal fans
coincide, in other words, they define the same toric variety (with
possibly different $\bQ$-polarizations). The following elementary lemma is
well-known, and we omit the proof.

\begin{lemma}\label{lem:toric-ex}
  \begin{enumerate}
  \item $P_{\w}$ and $P_{\w'}$ are normally equivalent iff $\w'$
    belongs to the interior $\cham(\w)$ of a certain rational
    polytope.
  \item If $\w'\in \overline{\cham(\w)}$ then there exists a natural
    morphism $\pi_{\w,\w'}:Y_{\w}\to Y_{\w'}$; it is birational if $\w'>\w$. 
  \item One has $\pi^*(L_{\w'})= L_{\w} + \sum (b_i'-b_i)B_i$. Thus,
    any positive multiple of this $\bQ$-line bundle that is integral,
    is semiample.
  \end{enumerate}
\end{lemma}

We note that this is precisely the kind of data that appears in
Theorem~\ref{thm:red-moduli}. 
Now consider a projective scheme $Z\subset\bP$ as in the general setup
\ref{setup}, $Z=\Proj \cA[\wZ]$. 

\begin{definition}
  For each weight $\wei$, let $\cA[\wZ\lw]\subset \cA[\wZ]$ be the
  subalgebra generated by the monomials whose characters lie in
  $\wt\Delta\lw$. We define $Z\lw:=\Proj \cA[\wZ\lw]$.
\end{definition}

\begin{theorem}\label{thm:properties-of-Zbeta}
  \begin{enumerate}
  \item Every $Z\lw$-polytope is the intersection of a $Z$-polytope
    with~$\Delta\lw$.
  \item There exists a rational map $Z\dashrightarrow Z\lw$. It is
    regular on the open subset of points of $Z$ whose moment polytopes
    intersect $\Delta\lw$.
  \item $Z$ and $Z\lw$ share an open subset $Z\lw^0$ whose points are
    the points with moment polytopes intersecting $\Int\Delta\lw$.
  \item There exists a chamber decomposition of the weight domain into
    finitely many interiors of polytopes, with the following properties:
    \begin{enumerate}
    \item If $\wei,\wei'$ belong to the same chamber, then $Z\lw=Z\lwp$.
    \item If $\wei'\in\overline{\cham(\wei)}$ then there exists a
      proper morphism $Z\lw\to Z\lwp$.
    \end{enumerate}
  \item Further, assume that for any $\cA$-field $k$, the corresponding
    variety $Z_k$ is integral, normal, and its monomials span the
    whole $\wt\Delta$. Then $Z\dashrightarrow Z\lw$ is a birational
    map, and in the
    previous statement one can take the $Z$-chamber decomposition
    defined in \ref{defn:Z-chamber-decomp}.
  \end{enumerate}
\end{theorem}
\begin{proof}
  We reduce the proof to the elementary case \ref{lem:toric-ex} of
  toric schemes, as follows. The homomorphism $\phi:\bZ^N\to \bZ^n$ in
  the setup~\ref{setup} gives a surjective map of polytopes $\sigma\to
  \Delta$, where $\sigma$ is a simplex with $N$ vertices. The preimage
  of $\Delta\lw$ is a certain subpolytope $\sigma\lw \subset \sigma$.
  Monomials of high enough degree $d$ generate the subalgebra
  $\cA[\wZ\lw]^{(d)}$, and this gives the embedding of $Z\lw$ into the
  toric scheme corresponding to the polytope $\sigma\lw$.

  Now the properties (1--4) are elementary for the ambient toric
  schemes, and hence they also hold for the subschemes $Z\lw$.

  To prove (5), note that that (2) and (3) together imply that
  $Z\dashrightarrow Z\lw$ is a birational map. Let $\wei,\wei'$
  belong to the same $Z$-chamber. Then on every geometric fiber we get a
  birational morphism $\varphi_k:(Z\lw)_k \to (Z\lwp)_k$. Since the
  $Z\lw$-polytopes and $Z\lwp$-polytopes are the same, 
  the $T$-orbits of $Z\lw,Z\lwp$ are in a dimension-preserving
  bijection, and so $\varphi_k$ is
  finite. Since $(Z\lwp)_k$ is normal, $\varphi_k$ is an isomorphism by
  the Main Zariski theorem. Since $Z\lw$ and $Z\lwp$ are free over
  $\cA$, and  $\varphi:Z\lw \to Z\lwp$ is an isomorphism fiberwise, it
  is an isomorphism.
\end{proof}

We now specialize to the case of grassmanians. Thus,
for every weight $\wei\in\cD(r,n)$ we get a projective scheme
$\gr\lw$, which we will call the \emph{weighted grassmannian}, and the
collection $\{\gr\lw\}$ satisfies the conclusions of
Theorem~\ref{thm:properties-of-Zbeta}, where the chamber decomposition
is the one defined in Definition~\ref{defn:chambers}.

\section{GIT theory of the universal family over the grassmanian}
\label{sec:git}

A key role in our definition of weighted stable hyperplane
arrangements will be played by the Geometric Invariant Theory of the
universal family $U\to \gr(r,n)$. Let us first review the basics
relevant to our case.

Let $Z\subset \bP$ be as in the setup~\ref{setup}. Then we have an
action of $T=\wT/\diag(\bG_m)$ on $Z$ and an action of $\wT$ on each
$\cO_Z(d)$, $d\in\bN$. The character group of $T$ is $\chi(T)=\{
(x_i)\in\bZ^n \mid \sum x_i=0 \}$.

A $T$-linearization of $\cO_Z(d)$ is an extension of the $T$-action
from $Z$ to $\cO_Z(d)$. It can be given by assigning to a monomial
$z^m=z_1^{m_1}\dots z_n^{m_n}$ of degree $\sum m_i = d$ an element
$\phi(m)\in \chi(T)$ so that $\phi(z^m)-\phi(z^{m'}) = m-m'$. This is
equivalent to choosing an element $\wei\in\bZ^n$ of degree $d$, so
that $\phi(z^m)=m-\wei$. 
We can also take $\wei$ to be of arbitrary positive degree, and
subtract the unique element of degree $d$ on the line $\bQ\wei$.
Then for every $d'$ that is a multiple of $d$, the
element of degree $d'$ on this ray also describes the induced linearization
of $\cO_Z(d')$.

Given a $T$-linearized ample sheaf $L=\cO_Z(d)$, one considers the ring of
sections $R=R(Z,L)=\oplus_{d\ge0} \Gamma(Z,L^d)$. This ring was
already graded by $\bZ^n$ by the setup. The linearization provides a
new grading by $\chi(T)=\bZ^{n-1}\subset \chi(\wT)=\bZ^n$. 

The GIT quotient $Z\lwq T$ is defined to be $\Proj R\lw$, where
the latter denotes the elements of degree 0 in the
$\chi(T)$-grading. In the original $\bZ^n$-grading, this means that we
consider the elements spanned by the monomials whose character in
$\bZ^n$ lies on the line $\bQ\wei$. Note as well that replacing $L$ by
a positive power does not change $Z\lwq T$. Hence, the input for this
construction is a weight $\wei$ up to a multiple, and an ample
invertible sheaf $L$ up to a multiple. 

Applied to the grassmannian $\gr(r,n)$ and the Pl\"ucker line bundle
$\cO_{\gr}(1)$, this means that every weight $\wei\in\cD(r,n)$ gives a
linearization and a GIT quotient $\gr\lwq T$. The quotients do not
respect the chamber structure of $\cD(r,n)$, however.

\emph{Our key observation now is that the chamber structure describes
  not the GIT quotients of $\gr$ but those of the universal family $U$
  over it.}

Let $U\subset \bP^{n-1}\times \gr(r,n)$ be the universal family of linear
spaces $\bP V \subset \bP^{n-1}$. Each of the $n$ hyperplanes
$H_i=\{z_i=0\}$ in $\bP^{n-1}$ defines a hyperplane $B_i \subset \bP V$,
unless $\bP V\subset H_i$.

The natural ample invertible sheaves on $U$ are 
$L_{a,b}=p_1^*\cO_{\bP^{n-1}}(a) \otimes p_2^*\cO_{\gr}(b)$ for
$a,b\in\bN$. The total ring of global sections of all of them is 
\begin{displaymath}
  \oplus_{a,b\ge0}H^0(U,L_{a,b}) = \frac{\cA[z_i,p_I]}
  {( \text{Pl\"ucker relations on }p_I,\ r_J )}
\end{displaymath}
where for each $J=\{j_0,\dotsc,j_r\}$, $r_J = \sum (-1)^k z_{j_k}
p_{J\setminus j_k} $. The $\bZ^n$-character of each $z_i$ is $e_i$, the
$i$-th coordinate vector in $\bZ^n$, and for $p_I$ it is $\sum_{i\in
  I}e_i$. Hence, the $\bZ$-degrees of $z_i$ and $p_I$ are $1$ and $r$
respectively, differing slightly from \ref{setup}.

\begin{definition}
  We choose:
\begin{enumerate}
\item The ample $\bQ$-line bundle $L_{a,b}$ with $(a,b)=(|\wei|-r,1)$,
  or any actual ample invertible sheaf for a multiple $(ma,mb)$ such
  that $m\wei$ is integral. 

  Note that if $\bP V\not\subset \cup H_i$ then $\cO(K_{\bP V}+ \sum
  b_i B_i) = \cO_{\bP V}(|\wei|-r)$.
\item The $T$-linearization corresponding to $\wei$.
\end{enumerate}
We denote the corresponding GIT quotient $U//\lw T$ by $\gr\lw^e$.
\end{definition}

\begin{lemma}\label{lem:gr-e}
  $\gr\lw^e$ is a closed subscheme of $\gr\lw$.
\end{lemma}
\begin{proof}
  Suppose that $m\wei$ is integral, and restrict to the subalgebra
  $\cA[\wgr\lw]^{(m)}$ whose homogeneous elements have degrees
  divisible by $m$. 
  A monomial in $p_I$ can be complemented to a monomial in $p_I,z_i$
  whose character is proportional to $\wei$ exactly when its
  character, divided by the number of $p_I$'s, lies in $\Delta(r,n)$ and
  in the cone $\wei-\bR^n_{\ge0}$. The intersection of these two sets
  is precisely $\Delta\lw$.

  Hence, we have a surjective homomorphism $\cA[\wgr\lw]^{(m)}\to
  R^{(m)}\lw$ sending a monomial $\prod p_I$ to its complement
  $z^s\prod p_I$. This gives the closed embedding.
\end{proof}

\begin{example}
If $\wei=\bo$ then every monomial in $p_I$ can be
complemented, and $\gr^e\subset\gr$ is the zero set of the equations
$r_J(p_I,\bo)$. Thus, $\gr^e_{\bo}$ is the same as $\gr^e$ that appeared in
Section~\ref{sec:unweighted-case}.
\end{example}

GIT gives the description of $U\lwq T$ in terms of orbits for each
geometric fiber. To recall,
there are two open subsets in $U$:
\begin{enumerate}
\item The set $U\uss\lw$ of semistable points $p$ for which there
  exists a section
  $$
  s\in R\lw = \oplus_{(a,b)\in\bQ\wei} H^0(U,L_{a,b})
  \text{ such that }  s(p)\ne 0.
  $$
\item The set $U\us\lw$ of (properly) stable points whose orbit in
  $U\uss\lw$ is closed and the stabilizer is finite; in our case
  trivial. (This set was denoted by $U\us_{(0)}$ in
  \cite{Mumford_GIT}. We use the currently prevalent notation.)
\end{enumerate}

Then we have a surjective morphism $U\lw\uss\to U\lwq T$,
the action is free on $U\us\lw$ and $U\lw\us/T$ is a geometric
quotient. Points of $U\uss\lw$ have the same image iff the
closures of their orbits intersect. Among such orbits, there
exists a unique closed one.

For the torus action, the Hilbert-Mumford's criterion for
(semi)stability takes an especially simple form. The following criterion
is well-known (e.g., cf. \cite{BrionProcesi}):
\begin{enumerate}
\item $p\in U\lw\uss \iff$ $\wei$ belongs to the moment polytope of $p$.
\item $p\in U\lw\us \iff$ $\wei$ lies in the interior of the moment
  polytope of $p$ and the latter is maximal-dimensional.
\end{enumerate}

The moment polytope here lies in $\bR^{n-1}$ which we shift so that it
lies in the hyperplane $\sum x_i = |\wei|$ in $\bR^n$. For our choice
$L_{|\wei|-r,1}$ of an ample $\bQ$-line bundle, the moment polytope of
the point $[p\in \bP V\subset \bP^{n-1}]\in U$ is:
\begin{displaymath}
  P_V + (|\wei|-r) \sigma_p,\quad \text{where}
\end{displaymath}
\begin{enumerate}
\item $P_V$, as before, is the matroid polytope of
  $[V\subset\bA^n]$, and
\item Denoting $I(p)=\{i \mid z_i(p)=0\}$, 
  $$
  \sigma_p = \left\{ (x_i)\in \bR^n \mid x_i\ge0, \ \sum x_i =1,
  \text{ and } x_i=0 \text{ for } i\in I(p) \right\}
  $$
\end{enumerate}

\begin{definition}
  $\Delta\lw^p$ is the face of $\Delta\lw$ where $x_i=b_i$ for $i\in I(p)$.
\end{definition}

\begin{lemma}
  \begin{enumerate}
  \item $p\in U\uss\lw \iff P_V \cap \Delta\lw^p \ne\emptyset.$
  \item $p\in U\us\lw \iff \Int P_V \cap \Int\Delta\lw^p \ne\emptyset$ and
    $P_V+\Delta\lw^p$ spans $\bR^{n-1}$.
  \end{enumerate}
\end{lemma}
\begin{proof}
  (1) $   \wei\in P_V + (|\wei|-r) \sigma_p \iff
  P_V \cap \{\wei - (|\wei|-r)\sigma_p\} \ne\emptyset $. The intersection
  of the latter polytope with $\Delta$ is $\Delta\lw^p$.
  
  (2) The point is stable iff we can replace $\wei$ with any nearby
  $\wei'$. This means that $\Int P_V \cap \Int\Delta\lw^p
  \ne\emptyset$ and $P_V+\Delta\lw^p$ spans $\bR^{n-1}$.
\end{proof}

\begin{theorem}\label{thm:git-lc-klt}
  \begin{enumerate}
  \item If $P_V\cap\Delta\lw=\emptyset$ or $V\subset \{z_i=0\}$ for
    some $i$ then no $p\in V$ is $\wei$-semistable.
  \item Suppose $P_V\cap\Delta\lw\ne\emptyset$. Then 
    $p\in U\uss\lw \iff (\bP V,\sum b_iB_i)$ is lc at $p$.
  \item Suppose $P_V\cap\Int\Delta\lw\ne\emptyset$. Then 
    $p\in U\us\lw \iff (\bP V,\sum b_iB_i)$ is klt at $p$.
  \end{enumerate}
\end{theorem}
\begin{proof}
  (1) If $P_V\cap\Delta\lw=\emptyset$ then $p\not\in U\uss\lw$ by the
  lemma. If $V\subset \{z_i=0\}$ then $P_V\subset\{x_i=0\}$, and
  $\Delta\lw^p \subset \{x_i=b_i\}$, so they do not
  intersect.

  (2) Suppose $P_V \cap \Delta\lw^p \ne\emptyset.$ Take $\alpha=(a_i)$
  in this intersection. By Theorem~\ref{thm:matroid-lc} the pair $(\bP
  V, \sum a_iB_i)$ is lc. Since one has $\sum a_i B_i = \sum b_i B_i$
  near $p$, the latter divisor is lc as well. 

  Vice versa, assume that $(\bP V, \sum b_iB_i)$ is lc at $p$. 
  By assumption, there exists 
  $\alpha\in P_V\cap\Delta\lw$. If $\alpha\not\in\Delta\lw^p$ then
  we are going to construct another $\alpha'=(a_i')\in
  P_V\cap\Delta\lw^p$. 

  If $P_V$ is maximal-dimensional and $\alpha\in\Int P_V$ then begin
  by increasing $x_i$ for $i\in I(p)$ until we get to $x_i=b_i$ while
  decreasing $x_i$ with $i\not\in I(p)$ and keeping $x_i\ge0$. This
  is possible to do since $\sum_{i\in I(p)} b_i \le \codim \cap_{i\in
    I(p)} B_i\le r-1$. By doing this, we either achieve the required
  $\alpha'$ or get to a lower-dimensional matroid polytope $P_{V'}$.
  But by Theorem~\ref{thm:matroid-polytopes}
  $P_{V'}$ is the product of maximal-dimensional polytopes for lower
  $(r_j,n_j)$. We finish by induction on $r$.

  (3) is proved by the same argument using the second part of 
  Theorem~\ref{thm:matroid-lc}. 
\end{proof}

\begin{definition}
  Denote by $\bP\lw$ the projective toric scheme over $\cA$ (a toric
  variety when working over $k$) corresponding to the polytope $\Delta\lw$.
\end{definition}

In particular, $\bP_{\bo}$ is the toric variety corresponding to the
hypersimplex $\Delta(r,n)$.

\begin{theorem}
  The morphism $U\uss\lw \to U\lwq T$ factors through $\bP\lw\times\gr\lw$.
\end{theorem}
\begin{proof}
  Consider $\bP^{n-1}\times\gr$ with the very ample sheaf $L_{a,b}$,
  $(a,b)\in\bQ\wei$. The rational map $\bP^{n-1}\times\gr
  \dashrightarrow \bP\lw \times \gr\lw$ is given by the monomials
  $\prod p_I z^m$ whose character is proportional to $\wei$ and, when
  normalized, the $\prod p_I$-part belongs to $\Delta\lw$ and the
  $z^m$-part belongs to $\beta-\Delta\lw$, which is just another copy
  of $\Delta\lw$, reflected. 

  This rational map is regular on the open subset where at least one
  of these monomials, considered as a section of $L^d$, is
  nonzero. But the ring generated by these monomials contains $R\lw$,
  so this open subset contains $U\uss\lw$.
\end{proof}

Recall from Section~\ref{sec:weight-grassm} that we denoted by
$\gr\lw^0$ an open subset of $\gr\lw$ and
$\gr$ that corresponds to $\Int\Delta\lw$. 

\begin{definition}
  $U\lw^0 \to \gr\lw^0$ will denote the pullback of $U\uss\to
  \gr\lw$ under the open inclusion $\gr\lw^0\to\gr$.
\end{definition}

\begin{theorem}\label{thm:free-git-quot}
  \begin{enumerate}
  \item The $T$-action on $U\lw^0$ is free.
  \item The geometric quotient $\hat\gr\lw^e = U\lw^0/T$ is projective
    and comes with a semiample invertible sheaf defining a proper
    birational morphism $\hat\gr\lw^e\to\gr\lw^e$, that is an
    isomorphism over $\gr\lw^e \cap \gr\lw^0$.
  \end{enumerate}
\end{theorem}
\begin{proof}
  (1) A point $V\in \gr\lw^0$ corresponds to a linear space
  $\bP V\subset\bP^{n-1}$ whose matroid polytope intersects
  $\Int\Delta\lw$. Take $p\in V$ such that $p\in U\uss\lw$. Then by
  Theorem~\ref{thm:git-lc-klt}(1) the pair $(\bP V,\sum b_iB_i)$ is lc
  at $p$. If we take a nearby $\wei'=(b_i')$ with $b_i'<b_i$ then
  $(\bP V,\sum b_iB_i)$ will be klt. Then by \ref{thm:git-lc-klt}(2)
  we have $p\in U\us\lwp$. Hence, $U^0\lw = U\us\lwp$, the action
  is free, and the quotient is projective.

  By removing $(U\uss\lw\setminus U^0\lw)$, we changed the
  equivalence relation on $U\lw^0$: for some of the orbits in $U\lw^0$
  their closures in $U\uss\lw$ intersect, and so they map to the same
  point of $U\lwq T$. The criterion of Theorem~\ref{thm:git-lc-klt}
  implies that every closed orbit in $U\uss\lw$ is contained in the
  closure of an
  orbit of $U\lw^0$. Hence, $\hat \gr\lw^e \to \gr\lw^e$ is
  surjective. It is given by the pullback of an ample invertible sheaf
  $\cO(m)$ on $\Proj R\lw$. This morphism is an isomorphism on the
  open subset $U\us\lw/T$ which contains $\gr\lw^e\cap \gr\lw^0$.
\end{proof}

\begin{remark}
  In the case of $\wei=\bo$ our construction is different from that of
  \cite{HackingKeelTevelev}. To explain it succinctly,
  \cite{HackingKeelTevelev} proceeds ``horizontally'', while we proceed
  ``vertically''. The points in $U\uss_{\bo}\setminus U\us_{\bo}$ are:
  the points $p\in \cup B_i$ 
  and the points $p\in V$ such that $P_V\subset \{x_i=1\}$
  for some $i$, satisfying the conditions of Theorem~\ref{thm:git-lc-klt}.
  So the action of $T$ on $U\uss_{\bo}$ is not
  free. There are several ways to restrict it to a subset with a free
  action:
  \begin{enumerate}
  \item ``Horizontally'', by removing the points $p\in \cup B_i$. The
    remaining set then is $U\cap(\bG_m^{n-1}\times\gr)$, where
    $\bG_m^{n-1}=\bP^{n-1}\setminus \cup H_i$. This is the choice of
    \cite{HackingKeelTevelev}. 
  \item ``Vertically'', by removing points with $P_V\subset \{x_i=1\}$
    -- our choice.
  \end{enumerate}
\end{remark}

\section{Definitions of the moduli space and the family}
\label{sec:defin-moduli-univ}

\begin{definition}
  (over $k=\bar k$)
  For each stable toric variety $Y\to \gr\lw$ over $\gr\lw$, we define
  the corresponding \emph{weighted stable hyperplane arrangement} as
  \begin{displaymath}
    X := (Y\times_{\gr\lw}U\uss\lw)\lwq T.
  \end{displaymath}
  We also define divisors $\bar B_i = (H_i\times\gr\lw) \cap U\uss\lw$ and
  then
  \begin{math}
    B_i := (Y\times_{\gr\lw} \bar B_i)\lwq T.
  \end{math}
\end{definition}

\begin{theorem}\label{thm:X-Xhat-Y}
  \begin{enumerate}
  \item The $T$-action on the restriction to $Y\times_{\gr\lw^0}U^0\lw$ is
    free. The geometric quotient $(\hat X, \hat B_i)$ 
    by this free action is projective and comes
    with a semiample invertible sheaf defining a proper birational morphism
    $\hat X\to X$, an isomorphism on the complements of $\hat B_i$, $B_i$.
  \item $X=Y\cap \gr\lw^e$.
  \item $X$ is reduced and $B_i$ are reduced Weil divisors on $X$.
  \end{enumerate}
\end{theorem}
\begin{proof}
  (1) follows immediately from Theorem~\ref{thm:free-git-quot} and (2) from
  Lemma~\ref{lem:gr-e}  by functoriality of GIT quotients. (3) follows
  since GIT quotients of reduced schemes are reduced.
\end{proof}

\begin{example}\label{exmp:XY-for-ha}
  Let $(\bP V,\sum b_iB_i)$ be an lc hyperplane arrangement. Then
  $\Delta\lw \subset P_V$. The weighted moment polytope of $(\bP
  V,\sum b_iB_i)$ is therefore $\Delta\lw$ itself. The
  normalization $Y$ of the closure of the orbit $T.V$ is a toric variety
  and it comes with a finite morphism to $\gr\lw$.
  
  The pullback of $\gr\lw^0$ under $Y\to \gr\lw$ is simply the orbit $T.V$,
  isomorphic to~$T$. Then $Y\times_{\gr\lw}U^0\lw = \bP V\times T$, and the
  quotient $\hat X$ is $\bP V$ itself, together with the original divisors
  $B_i$. Since $\hat X\simeq\bP^{r-1}$,  the morphism $\hat X\to X$
  must be an isomorphism. Hence, every lc hyperplane arrangement
  appears as a particular case of our construction.
\end{example}

\begin{theorem}
  \label{thm:wsha-sings}
  \begin{enumerate}
  \item   Any weighted stable hyperplane arrangement $(X,\sum b_i B_i)$ is a
  stable pair, i.e. it has slc singularities and $K_X+\sum b_iB_i$ is
  an ample $\bQ$-divisor.
\item $\hat X$   is Gorenstein.
  \item $X$ is Cohen-Macaulay, and $X\setminus \cup B_i$ is Gorenstein.
  \item Assume $mb_i\in\bZ$. Then 
    $m(K_X+\sum b_iB_i)$ is the
  restriction under $Y\subset X$ of the ample invertible sheaf $F_Y$ on $Y$
  corresponding to the polytope $m\Delta\lw$. 
  \end{enumerate}
\end{theorem}
\begin{proof}
(1)  By Theorems~\ref{thm:matroid-lc} and \ref{thm:git-lc-klt},
  $(U\lw^0,\sum b_i \bar B^0_i)$ is a family of open lc subsets of
  hyperplane arrangements. Hence, $U^0$ is smooth, and there exists a finite
  sequence of blowups of $\bP^{n-1}$ giving a simultaneous resolution
  of singularities of $(U\lw^0,\sum b_i \wB^0_i)$.

    On the other hand, a stable toric variety $Y$ together with its
    boundary is slc by \cite{Alexeev_CMAV}; and the boundary is
    contained in $\gr\lw\setminus\gr\lw^0$ so can be omitted.  
    The stable toric variety is Cohen-Macaulay, and its interior is
    Gorenstein by the Stanley-Reisner theory because the topological
    space $\Delta\lw$ is a smooth manifold with boundary.
    Therefore, the pullback $V:=Y\times_{\gr\lw}U^0\lw$, together with
    the boundary, has slc singularities, and it is Gorenstein. 

    Then the geometric quotient $\hat X=V/T$ by the free $T$-action is
    Gorenstein, giving (2).

    Now let $m\in\bN$ be such that $m\wei$ is integral, and let $F_V$,
    $F_{\hat X}$, $F_X$ be the invertible sheaves on $V$, $\hat X$,
    $X$ given by the GIT construction: $F_V$ is the pullback of
    $L_{m\wei-mr,m}$, sections of $F_V$ descend to
    sections of $F_{\hat X}$ and $F_X$, $F_{\hat X}$ is semiample and
    defines the contraction $\hat X\to X$, $F_X$ is ample. 

    We observe that by construction one has $F_{\hat X}=\cO_{\hat
      X}(m(K_{\hat X}+\sum b_i\hat B_i))$. This implies that
    $F_{X}=\cO_{X}(m(K_{X}+\sum b_iB_i))$ and that $(X,\sum b_iB_i)$
    is slc. 

    Since $X\setminus\cup B_i = \hat X\setminus\cup \hat B_i$,
    $X_i\setminus \cup B_i$ is Gorenstein. $X$ is Cohen-Macaulay
    because it is the result of a log crepant contraction isomorphic
    outside of $\cup B_i$ and there exists a positive combination of
    $B_i$ which is Cartier. 

    Let $F_Y$ be the (integral) ample invertible sheaf corresponding
    to the polytope $m\Delta\lw$. Then by the same argument as in
    Theorem~\ref{lem:gr-e} sections of $F_Y$ restrict to sections of
    $F_X$. This gives (3).
\end{proof}

Let $Y\to\gr\lw$ be a stable toric variety over $\gr\lw$, and
$\uP\lw=\{P_{V,\wei}\}$ be the corresponding cover of $\Int\Delta\lw$. 
Each of these polytopes has a parent, so that $P_{V,\wei}=P_V\cap
\Delta\lw$. 

We denote by $Y[P_{V,\wei}]$ the corresponding projective toric
variety. We also denote by $\sigma_n$ the simplex $\{(x_i)\in \bR^n
\mid x_i\ge0, \ \sum x_i=r\}$. The corresponding to it toric variety is
$\bP^{n-1}$. If $\pwei$ is maximal-dimensional then toric geometry
gives a natural birational map
$\psi[\pwei]:\bP^{n-1}\dashrightarrow Y[P_{V,\wei}]$, an isomorphism
on the torus $\bG_m^{n-1}$. 

Now let $\pwei$ be a weighted matroid polytope of codimension $c$, and
let $P_V$ be its parent, a matroid polytope. Recall from
Theorem~\ref{thm:matroid-polytopes} that $P_V= \prod P_j$, the product
of maximal-dimensional polytopes for a subdivision $\{1,\dotsc,n\}=
\sqcup_{j=0}^c I_j$. Then we have a natural rational map
$\bP^{n-1}\dashrightarrow \prod \bP^{|I_j|-1}$ which on an open subset
is the quotient by a free action by $\bG_m^{c}$.

In this case, we denote by $\psi[\pwei]:\bP^{n-1}\dashrightarrow
Y[P_{V,\wei}]$ the latter rational map followed by the birational map
$\prod \bP^{|I_j|-1} \to Y[\pwei]$ corresponding to polytopes $\prod
\sigma_{|I_j|}$ and $\pwei$.

\begin{theorem}\label{thm:wsha-detailed-desc}
  The stratification of $X$ into irreducible components and their
  intersections is in a dimension-preserving,
  with a shift by $n-r$,  bijection with the tiling
  $\cup P_{V,\wei}$ of $\Int\Delta\lw$. Moreover, the closure of the
  stratum corresponding to $\pwei$ is the closure of the image of $\bP
  V$ under $\psi[\pwei]:\bP^{n-1}\dashrightarrow Y[\pwei]$. This
  rational map is regular on the open subset of $\bP V$ where $(\bP V,
  \sum b_iB_i)$ is lc. The image of this regular set gives a locally
  closed stratum in $X$ corresponding to $\pwei$.
\end{theorem}
\begin{proof}
  Each of the polytopes corresponds to an arrangement $\bP V\subset
  \bP^{n-1}$ so that the moment polytope $P_V$ intersects
  $\Int\Delta\lw$. Then we simply follow the $T$-orbit
  of $\bP V$ in $\gr\lw^0$, to the pullback in $U^0\lw$, to the quotient
  in $X$, then back to the irreducible component of $Y[\pwei]$ under
  the inclusion $X\subset Y$. 

  When $\pwei$ is maximal-dimensional, the orbit $T.V$ in $\gr\lw^0$
  has trivial stabilizer. Hence, under the quotient by the free action
  by $T$, the open subset of lc points on $\bP_V$ is preserved; then
  part is contracted by a birational morphism. If $\bP_V$ has codimension
  $c$ then the stabilizer of $T.V$ in $\gr\lw^0$ is $\bG_m^{c}$. Then
  the quotient factors through the quotient of an open subset of
  $\bP^{n-1}$ by $\bG_m^c$.
\end{proof}

We are now ready to define the moduli space and the universal family
of pairs over it. 

\begin{definition}\label{defn:my-moduli}
  $\oM\lw(r,n) = \M^{T}(\gr\lw,\Delta\lw).$ 
\end{definition}

\begin{lemma}
  $\M^{T}(\gr\lw,\Delta\lw)$ is a fine moduli space.
\end{lemma}
\begin{proof}
  Indeed, every $T$-orbit in $\gr\lw^0$ is also an orbit in $\gr$. If
  its stabilizer is finite then it is in fact trivial. Therefore,
  every irreducible component of a stable toric variety $Y\to \gr\lw$
  maps to its image birationally, and the automorphism group of $Y\to
  \gr\lw$ is trivial. Then $M^T$ is a fine moduli space.
\end{proof}

\begin{definition}\label{defn:my-family}
  The family of
  weighted stable hyperplanes arrangements $(\cX,\cB_i)\to
  \oM\lw(r,n)$ is the GIT quotient of the pullback of $U\lw^0\to
  \gr\lw$ by the universal family of stable toric varieties $\cY\to
  \M^{T}(\gr\lw,\Delta\lw).$
\end{definition}

\begin{theorem}
  $(\cX,\cB_i)\to\oM\lw(r,n)$ is a locally free (in particular, flat)
  morphism.   
\end{theorem}
\begin{proof}
  The families $U\to \gr$ and $\cY\to \oM\lw(r,n)$ are locally free,
  i.e. locally they are given by locally free modules. This implies
  that the pullback $\cY \times_{\oM} U\lw\uss$ is locally free.
  Algebraically, the GIT quotient is constructed by taking the
  degree-0 component in an algebra. Thus, this subalgebra is a direct
  summand, and a direct summand of a locally free module is locally
  free (by Kaplansky's theorem \cite{Kaplansky_Projective}, over any
  ring a module is locally free iff it is projective).
\end{proof}

\section{Completing the proofs of main theorems}
\label{sec:completing-proofs}

\begin{proof}[Proof of Theorem \ref{thm:existence} (Existence)]
  The parts (1) and (3) were established in the previous section. The
  subset $\M\lw(r,n) \subset \oM\lw(r,n)$ is the open subset of
  $\M^{T}(\gr\lw,\Delta\lw)$ where the stable toric varieties are
  irreducible, cf. Example~\ref{exmp:XY-for-ha}. The sheaf
  \linebreak
  $\cO_{\cX}(m(K_{\cX} + \sum b_iB_i))$ is free over $\oM\lw$ because
  by Theorem~\ref{thm:wsha-sings} it is the restriction of the
  invertible ample sheaf $F_{\cY}$ from the universal family of
  stable toric varieties that corresponds to the lattice polytope
  $m\Delta\lw$, and $F_{\cY}$ is free: it's sections give the finite
  morphism to $\gr\lw$.

  The remaining part (2) is
  proved in the Reconstruction Theorem~\ref{thm:reconstruction} below.
\end{proof}

\begin{proof}[Proof of Theorem \ref{thm:red-moduli} (Reduction morphisms)]
(1) For $\wei,\wei'$ in the same chamber, we have $\gr\lw=\gr\lwp$ by
  Theorem~\ref{thm:properties-of-Zbeta}, applied to
  grassmannians. Also the conditions for GIT (semi)stability are the
  same. So the moduli and the families are the same.

(2) If $\wei'\in \overline{\cham(\wei)}$, we have a reduction morphism 
  $\gr\lw\to\gr\lwp$ again by
  Theorem~\ref{thm:properties-of-Zbeta}. The third application of the
  same theorem gives the reduction morphism between the stable toric
  varieties $Y\lw,Y\lwp$ over $\gr\lw,\gr\lwp$. Finally, this gives in a
  canonical way the reduction morphisms between the pullbacks of the
  universal families and their GIT quotients.

  Each $\pi_{\wei,\wei'}$ is log crepant. That is because the
  morphism on the ambient stable toric varieties is given by pullback
  of $L_{\wei}+\sum(b_i'-b_i)B_i$ (cf.  \ref{lem:toric-ex},
  \ref{thm:properties-of-Zbeta}), and $L_{\wei}$ on $Y$ restricts to
  $K_X+\sum b_i B_i$ on $X$ by Theorem~\ref{thm:wsha-sings}.

(3)  When specializing up, the morphism $\oM\lw\to \oM\lwp$ is an
  isomorphism. Indeed, a stable toric variety $Y\to\gr\lw$ is uniquely
  determined by its restriction $Y^0$ to $\gr\lw^0$: $Y$ is the partial
  normalization at the boundary of the closure of $Y^0$. But
  $\gr\lw^0=\gr\lwp^0$ in this case. 

  Additionally, when specializing up, the morphism $X\lw\to X\lwp$ is
  simply our morphism $\hat X\to X$, so by Theorem~\ref{thm:X-Xhat-Y}
  it is birational and an isomorphism outside $\hat B_i,B_i$.

  (4) is an immediate consequence of the parts (1,2,3).
\end{proof}

Note that if the source of the GIT quotient were fixed, with only the
line bundle and the polarization changing, the statement would be an
application of the well-known theory of variation of GIT
quotients \cite{BrionProcesi, DolgachevHu}.

\begin{proof}[Proof of Theorem~\ref{thm:small-weights} (Moduli for
  small weights)] 
  In this case $\gr\lw^0=\gr\uss\lw=\gr\us\lw$, the $T$-action on it is
  free, and every stable toric variety over $\gr\lw$ is the closure of
  a unique $T$-orbit. Hence, $M^T(\gr\lw,\Delta\lw)=\gr\us\lw/T
  = \gr\lwq T$.

  The equivalence of the GIT quotients 
  $(\bP^{r-1})^n//\PGL(r)$ and $\gr(r,n)//T$ is well-known, see,
  e.g., \cite[2.4.7]{Kapranov_QuotientsGrassmannians}.
\end{proof}

  Intuitively, the contribution $p_1^* \cO_{\bP^{n-1}}(|\wei|-r)$ to 
  the polarization $L_{|\wei|-r,1}$ in
  this case approaches zero and only the quotient of $\gr\lw$ remains.

\begin{theorem}[Reconstruction Theorem]
  \label{thm:reconstruction}
  The stable toric variety $Y\to\gr\lw$ can be uniquely reconstructed
  from $(X,\sum b_iB_i)$.
\end{theorem}
\begin{proof}
  Let $Y[\pwei]$ be an irreducible component of a stable toric variety
  $Y\to \gr\lw$, as in Theorem~\ref{thm:wsha-detailed-desc}, and
  $X[\pwei] \subset Y[\pwei]$ be the corresponding irreducible
  component of $X$. We first show that $Y[\pwei]$ can be reconstructed from
  $X[\pwei]$ intrinsically.

  Indeed, the boundary of $X[\pwei]$ in $X$ is labelled by the
  divisors $B_i$, some of them coinciding. Then
  the defining inequalities of $\pwei$ can be read off this 
  configuration: every missing intersection $\cap_{i\in I} B_i$ of
  codimension $k$ gives the inequality $\sum_{i\in I}x_i\le k$. 
  This recovers the polytope $\pwei$. 

  Then the embedding $X[\pwei]\to Y[\pwei]$ is recovered as
  follows. For every $m$ such that $m\wei$ is integral, every integral
  point $(x_i)\in m\pwei$ gives a section of
  the sheaf $\cO_{X[\pwei]}( m(K_X+\sum b_iB_i) )$.
   Namely, it is a unique up to a constant section vanishing
  at $B_i$ to the order $x_i$. The collection of these sections gives
  the embedding $X[\pwei]\to Y[\pwei]$ defined up to $n$ choices of
  multiplicative constants, one for each $B_i$, i.e. up to the action
  of $\wT$.  

  Finally, $\bP V$ is recovered from the image of $X[\pwei]\to
  Y[\pwei]$ by applying the inverse of the rational map $\psi[\pwei]$
  of Theorem~\ref{thm:wsha-detailed-desc}. Then the orbit $T.V$ in
  $\gr\lw^0$ gives the morphism $Y[\pwei]\to \gr\lw$. The whole
  stable toric variety  $Y\to \gr\lw$ is recovered this way by looking
  at all maximal-dimensional polytopes $\pwei$.
\end{proof}

We note that for $\wei=\bo$ this proof is very different from the one
given in \cite{HackingKeelTevelev}, which does not extend to the
weighted case.

\section{Some simple examples}
\label{sec:some-examples}

\begin{example}
  $(r,n)=(2,4)$, $\wei=\bo$. Consider the subdivision of $\Delta(2,4)$,
  the octahedron on the ${4 \choose 2}$ vertices $ij$,
  into two pyramids: $P_1$ missing the vertex 34, and $P_2$, missing
  the vertex $12$.

  $P_1$ corresponds to the configuration of 4 points in $\bP^1$ for
  which the Pl\"ucker coordinate
  $p_{34}=0$, i.e. $B_3=B_4$. This polytope is given by the
  inequality $x_3+x_4\le 1$, which is precisely the lc condition for
  this configuration. 
  Similarly, for $P_2$ one has
  $B_1=B_2$. On the intersection $P_1\cap P_2$ one has $B_1=B_2$,
  $B_3=B_4$, and the defining inequalities become $x_1+x_2=1$,
  $x_3+x_4=1$, i.e. $P_1\cap P_2 = \Delta(1,2) \times \Delta(1,2)$. 

  The irreducible component $X_1$ then is the closure of the image of
  $\bP^1\setminus B_3$, so isomorphic to $\bP^1$; and similarly for
  $X_2$. The intersection $X_1\cap X_2$ is the quotient
  $(\bP^1\setminus \{B_1,B_3\})/\bG_m$, so a point. So $X$ is a union
  of two $\bP^1$'s intersecting at a point.
\end{example}

\begin{example}
  $(r,n)=(2,4)$, $\wei=(1/2,1/2,1,1)$. Consider the trivial
  subdivision of $\Delta\lw$, with just the polytope itself.
  The points $B_1$ and $B_2$ may or may not coincide, 
  depending on whether the parent polytope
  is the pyramid $P_2$ from the previous example, or $\Delta(2,4)$,
  otherwise the points are 
  pairwise distinct. $X=\bP V=\bP^1$. 

  Say, $B_1=B_2$. Then $X\subset Y$ intersects the stratum corresponding
  to the edge $x_1=x_2=1/2$ of $\Delta\lw$, at a point $q$. In this case, the
  $T$-translates of $X$ do not sweep out an open subset of $Y$, and this is
  very different from the unweighted situation of
  Section~\ref{sec:unweighted-case}. 
  
  If we consider the GIT quotient of the pullback family over the
  whole $Y$ (not just $Y\cap\gr\lw^0$ as in our construction), then on
  the boundary some fibers to the GIT quotient are modelled on the
  curve $\bA^1\cup_q \bA^1$, which is a transversal slice of $Y$
  at the point $q$.
\end{example}

In all cases with  $r=2$ the considerations are quite similar, and
produce a tree of $\bP^1$'s. 

\begin{example}
  $(r,n)=(3,5)$, this will correspond to
  Example~\ref{exmp:3-degs}. Begin with $\wei=\bo$, and consider the
  subdivion of $\Delta(3,5)$ into 3 polytopes:
  $P_0=\{x_1+x_2+x_5\le2, \ x_3+x_4+x_5\le2\}$, $P_1=\{x_1+x_2\le1\}$,
  and $P_2=\{x_3+x_4\le 1\}$.

  Then $P_0$ corresponds to the configuration of $5$ lines such that
  $B_1\cap B_2\cap B_5$ is a point, $B_3\cap B_4\cap B_5$ is a point,
  and otherwise generic. The matroid polytope $P_V$ is obtained from
  $\Delta(3,5)$ by cutting two corners, and the intersection
  $P_V\cap \{x_5=~1\}$ has codimension 2, not 1 as might be expected: it is 
  $\{x_5=1,x_1+x_2\le 1,x_3+x_4\le 1\}$, so the corresponding face gets
  contracted. 

  As in Theorem~\ref{thm:wsha-detailed-desc}, the irreducible
  component $X^0$ is the image of $\bP^2$ under $\bP^4\dashrightarrow
  Y[\pwei]$. On $\bP^2$ it blows up two points and contracts the
  strict preimage of~$B_5$. 
  The configuration $(\bP V,\sum B_i)$ is lc outside of two points, so
  the divisor $B_5$ is present on $\hat X$; it is contracted by the
  log crepant morphism $\hat X\to X$. 

  For the weight $\wei=(1,1,1,1,1-\epsilon)$, the face $P_V\cap
  \{x_5=1-\epsilon \}$ has codimension 1, and the curve $B_5$ is not
  contracted. 
\end{example}

\begin{example}
  Consider the subdivision of $\Delta\lw$ 
  by a single hyperplane $x_1
  + \dotsb + x_{n_1} = r_1$, equivalently $x_{n_1+1}+\dotsb +
  x_n=r_2$, with $r_1+r_2=r$, $n_1+n_2=n$. Then $X$ is the union of
  $\Bl_{\bP^{r_1-1}}\bP^{r-1}$  and $\Bl_{\bP^{r_2-1}}\bP^{r-1}$ glued
  along   $\bP^{r_1-1} \times \bP^{r_2-1}$. 
\end{example}

\begin{example}
  Let $(a_1,\dotsc,a_{n-r+1})\in\cD(1,n-r+1)$ be a weight such that
  $\sum a_i>1$ but $\sum_{i\in I} a_i\le 1$ for any proper subset
  $I$. Let $\wei\in\cD(r,n)$ be the weight consisting of $\alpha$
  preceded by $(r-1)$ $1$'s.

  Then $\M\lw=\oM\lw= (\bP^{n-r-1})^{r-1}$. For $r=2$ this was
  established in \cite[4.5]{AlexeevGuy}. For the general case, we
  first observe that Theorem~\ref{thm:small-weights} applies in this
  case after replacing $1$ with $1-\epsilon$ for some $0<\epsilon\ll1$,
  and so $\oM\lw$ a moduli space of lc hyperplane
  arrangements. The lc condition implies that the $(r-1)$ hyperplanes
  with weight $1$ must intersect normally. Restricting to an
  intersection to any $(r-2)$ of these hyperplanes, a line, gives the $r=2$
  situation, for the weight $(1,a_i)$, and the moduli space for this
  is $\bP^{n-r-1}$. Each of the hyperplanes with weight $a_i$ is
  uniquely determined by the intersections with these $(r-1)$ lines,
  and all of these configurations are lc. So $\M\lw= (\bP^{n-r-1})^{r-1}$. 
\end{example}

\begin{example}
  Let $\wei=(1,\dotsc,1,\epsilon,\dotsc,\epsilon)$,
  $|\wei|=r+(n-r)\epsilon$. The case of $r=2$ was introduced in
  \cite{LosevManin}, and $\oM\lw(2,n)$ is the toric variety for the
  permutohedron, see also \cite[2.11(4)]{AlexeevGuy}.

  For any $r$, the closure of $\M\lw$ in $\oM\lw$ is the toric variety
  for the fiber polytope $\Sigma(\sigma_r^{\oplus (n-r)} \to
  (n-r)\sigma_r)$, where $\sigma_r$ is the simplex with $r$ vertices
  and side 1, and $(n-r)\sigma_r$ is $\sigma_r$ dilated by $(n-r)$.

  This moduli space also has an interpretation as the moduli space of
  stable toric pairs $(X,D_1,\dotsc,D_{n-r})$, as in
  \cite{Alexeev_CMAV} but with $(n-r)$ divisors instead of
  one. Explaining this in detail would take quite some space, and is
  better done elsewhere.
\end{example}

\bibliographystyle{amsalpha}

\begin{thebibliography}{GGMS87}

\bibitem[AB04a]{AlexeevBrion_Affine}
V.~Alexeev and M.~Brion, \emph{Stable reductive varieties. {I}. {A}ffine
  varieties}, Invent. Math. \textbf{157} (2004), no.~2, 227--274.

\bibitem[AB04b]{AlexeevBrion_Projective}
\bysame, \emph{Stable reductive varieties. {II}. {P}rojective case}, Adv. Math.
  \textbf{184} (2004), no.~2, 380--408.

\bibitem[AB06]{AlexeevBrion_SphericalModuli}
\bysame, \emph{Stable spherical varieties and their moduli}, IMRP Int. Math.
  Res. Pap. (2006), Art. ID 46293, 57.

\bibitem[AG06]{AlexeevGuy}
V.~Alexeev and G.~M. Guy, \emph{Moduli of weighted stable maps and their
  gravitational descendants}, Journal of the Institute of Mathematics of
  Jussieu, to appear, arXiv: math.AG./0607683.

\bibitem[AK06]{AlexeevKnutson}
V.~Alexeev and A.~Knutson, \emph{Complete moduli spaces of branchvarieties},
  arXiv math.AG/0602626 (2006).

\bibitem[Ale02]{Alexeev_CMAV}
V.~Alexeev, \emph{Complete moduli in the presence of semiabelian group action},
  Ann. of Math. (2) \textbf{155} (2002), no.~3, 611--708.

\bibitem[Ale06]{Alexeev_ICMtalk}
\bysame, \emph{Higher-dimensional analogues of stable curves}, Proceedings of
  the International Congress of Mathematicians, Vol. II (Madrid, 2006),
  European Math. Soc. Pub. House, 2006.

\bibitem[BP90]{BrionProcesi}
M.~Brion and C.~Procesi, \emph{Action d'un tore dans une vari\'et\'e
  projective}, Operator algebras, unitary representations, enveloping algebras,
  and invariant theory (Paris, 1989), Progr. Math., vol.~92, Birkh\"auser
  Boston, Boston, MA, 1990, pp.~509--539.

\bibitem[DH98]{DolgachevHu}
I.~V. Dolgachev and Y.~Hu, \emph{Variation of geometric invariant theory
  quotients}, Inst. Hautes \'Etudes Sci. Publ. Math. (1998), no.~87, 5--56,
  With an appendix by Nicolas Ressayre.

\bibitem[GGMS87]{GelfandGoreskyMacPhersonSerganova}
I.~M. Gel{\cprime}fand, R.~M. Goresky, R.~D. MacPherson, and V.~V. Serganova,
  \emph{Combinatorial geometries, convex polyhedra, and {S}chubert cells}, Adv.
  in Math. \textbf{63} (1987), no.~3, 301--316.

\bibitem[Has03]{Hassett_WeightedStableCurves}
B.~Hassett, \emph{Moduli spaces of weighted pointed stable curves}, Adv. Math.
  \textbf{173} (2003), no.~2, 316--352.

\bibitem[HKT06]{HackingKeelTevelev}
P.~Hacking, S.~Keel, and J.~Tevelev, \emph{Compactification of the moduli space
  of hyperplane arrangements}, J. Algebraic Geom. \textbf{15} (2006), no.~4,
  657--680. \MR{MR2237265 (2007j:14016)}

\bibitem[HS04]{HaimanSturmfels}
M.~Haiman and B.~Sturmfels, \emph{Multigraded {H}ilbert schemes}, J. Algebraic
  Geom. \textbf{13} (2004), no.~4, 725--769.

\bibitem[Kap58]{Kaplansky_Projective}
I.~Kaplansky, \emph{Projective modules}, Ann. of Math (2) \textbf{68} (1958),
  372--377.

\bibitem[Kap93]{Kapranov_QuotientsGrassmannians}
M.~M. Kapranov, \emph{Chow quotients of {G}rassmannians. {I}}, I. M.
  Gel{$'$}fand Seminar, Adv. Soviet Math., vol.~16, Amer. Math. Soc.,
  Providence, RI, 1993, pp.~29--110.

\bibitem[Laf03]{Lafforgue_Chirurgie}
L.~Lafforgue, \emph{Chirurgie des grassmanniennes}, CRM Monograph Series,
  vol.~19, American Mathematical Society, Providence, RI, 2003.

\bibitem[LM00]{LosevManin}
A.~Losev and Y.~Manin, \emph{New moduli spaces of pointed curves and pencils of
  flat connections}, Michigan Math. J. \textbf{48} (2000), 443--472, Dedicated
  to William Fulton on the occasion of his 60th birthday.

\bibitem[MFK94]{Mumford_GIT}
D.~Mumford, J.~Fogarty, and F.~Kirwan, \emph{Geometric invariant theory}, third
  ed., Ergebnisse der Mathematik und ihrer Grenzgebiete (2) [Results in
  Mathematics and Related Areas (2)], vol.~34, Springer-Verlag, Berlin, 1994.

\bibitem[PS02]{PeevaStillman}
I.~Peeva and M.~Stillman, \emph{Toric {H}ilbert schemes}, Duke Math. J.
  \textbf{111} (2002), no.~3, 419--449.

\end{thebibliography}

\def\cprime{$'$} \def\cprime{$'$}
\providecommand{\bysame}{\leavevmode\hbox to3em{\hrulefill}\thinspace}
\providecommand{\MR}{\relax\ifhmode\unskip\space\fi MR }
\providecommand{\MRhref}[2]{%
  \href{http://www.ams.org/mathscinet-getitem?mr=#1}{#2}
}
\providecommand{\href}[2]{#2}

\end{document}